\documentclass[10pt]{amsart}
\pdfoutput=1
\usepackage[utf8]{inputenc}
\usepackage[T1]{fontenc}
\usepackage[full]{textcomp}

\usepackage{csquotes}
\usepackage[english]{babel}

\usepackage{mathtools}
\mathtoolsset{mathic} %

\usepackage{cfr-lm}
\usepackage{dsfont}

\usepackage{microtype}

\usepackage{amsmath}
\usepackage{amssymb}
\usepackage{amsthm}

\usepackage[inline]{enumitem}

\usepackage[usenames,dvipsnames]{xcolor}

\usepackage{tikz-cd}

\usepackage[backend=biber,doi=false,url=false,isbn=false,%
safeinputenc,sorting=nyt]{biblatex}
\bibliography{\jobname.bib}
\bibliography{extra.bib}

\newtheorem*{theorem-non}{Theorem}

\theoremstyle{definition}

\newtheorem*{conjectures-non}{Conjectures}
\newtheorem*{defition-non}{Definition}
\newtheorem*{remark-non}{Remark}

\def\into{\hookrightarrow}

\def\id{\textnormal{id}}

\def\NN{\mathbb{N}}
\def\ZZ{\mathbb{Z}}
\def\QQ{\mathbb{Q}}

\def\CC{\mathbb{C}}

\def\QQl{\mathbb{Q}_\ell}

\def\Hom{\textnormal{Hom}}
\def\Aut{\textnormal{Aut}}

\def\Gal{\textnormal{Gal}}

\def\Vect#1{\textnormal{Vect}_{#1}}
\def\Rep#1{\textnormal{Rep}_{#1}}
\def\HS#1{\textnormal{HS}_{#1}}
\def\Mot#1{\textnormal{Mot}_{#1}}

\def \HH{\textnormal{H}}

\newif\ifcc\def\cc{\cctrue}
\newif\ifsp\def\sp{\sptrue}

\def\group#1#2{#1_{\ifsp{1,}\fi#2}^{\ifcc\circ\fi}\ccfalse\spfalse}

\def\mG#1{\group{\mathcal{G}}{#1}}    %
\def\rG#1{\group{\textnormal{G}}{#1}} %

\def\rGs{\rG\sigma}
\def\rGl{\rG{\ell,k}}
\def\rGlc{\cc\rGl}

\def\motTC{\textnormal{TC}'}
\def\MST{\textnormal{MST}}
\def\MTC{\textnormal{MTC}}
\def\AST{\textnormal{AST}}

\def\GL{\textnormal{GL}}
\def\Gm{\mathbb{G}_{\textnormal{m}}}

\def\onto{\twoheadrightarrow}
\def\Im{\textnormal{Im}}

\begin{document}
\title[Mumford--Tate implies algebraic Sato--Tate]
{The Mumford--Tate conjecture implies\\
the algebraic Sato--Tate conjecture\\
of Banaszak and Kedlaya}
\author[V. Cantoral-Farf\'an]{
Victoria Cantoral-Farf\'an
}
\address[Victoria Cantoral-Farf\'an]{
\begin{itemize}
\item[-] KU Leuven, Department of Mathematics\\
Celestijnenlaan 200B, B-3001 Leuven, Belgium
\end{itemize}
} \email[Victoria Cantoral-Farf\'an]{victoria.cantoralfarfan@kuleuven.be}

\author[J. Commelin]{
Johan Commelin
}
\address[Johan Commelin]{
\begin{itemize}
\item[-] Albert-Ludwigs-Universität Freiburg, Mathematisches Institut\\
 Ernst-Zermelo-Stra\ss{}e 1, 79104 Freiburg im Breisgau, Deutschland
\end{itemize}
} \email[Johan Commelin]{jmc@math.uni-freiburg.de}

\date{\today}
\thanks{The first author was partially supported by KU Leuven IF C14/17/083. The second author was supported by the
Deutsche Forschungs Gemeinschaft~(DFG)
under Graduiertenkolleg~1821 (\emph{Cohomological Methods in Geometry}).}
\subjclass[2010]{Primary \textsc{14f20}; Secondary \textsc{11g10}}
\keywords{Algebraic Sato--Tate conjecture,
 Mumford--Tate conjecture, abelian varieties, motives}

\begin{abstract}
 The algebraic Sato--Tate conjecture was initially introduced by Serre
 and then discussed by Banaszak and Kedlaya
 in the papers~\cite{MR3320526} and~\cite{MR3502937}.
 This note shows that the Mumford--Tate conjecture for an abelian variety~$A$
 implies the algebraic Sato--Tate conjecture for~$A$.

 The relevance of this result lies mainly in the fact
 that the list of known cases of the Mumford--Tate conjecture
 was up to now a lot longer than the list of known cases of the
 algebraic Sato--Tate conjecture.
\end{abstract}

\maketitle

\section{Introduction}

In 1966, Serre wrote a letter~\cite{SerreLettre} to Borel
where he presented some remarkable links
between the Mumford--Tate conjecture and questions
related to the equidistribution of traces of Frobenius.
Inspired by the previous work of Serre,
the algebraic Sato--Tate conjecture was introduced by Banaszak and Kedlaya
as an attempt to prove new instances
of the generalized Sato--Tate conjecture~\cite{MR3320526,MR3502937}.
These articles are considered as the theoretical motivation that allowed
Fit\'e, Kedlaya, Rotger and Sutherland to provide the classification
of all the possible Sato--Tate groups that can appear
for the Jacobian of a curve of genus~$2$ over a number field~\cite{FKRS}.

For more details about the generalized Sato--Tate conjecture
and the conjectural relation with the Mumford--Tate group
we refer to the presentations of~\cite{Fite,MandM},
\S13 of Serre's paper~\cite{Serre94}, and~\cite{SerreNp}.

\medskip\noindent
\textbf{Main result.}\quad
In this paper we show that for abelian varieties
(in fact, abelian motives)
the Mumford--Tate conjecture implies the algebraic Sato--Tate conjecture.

In the next section we recall definitions,
introduce some notation, and derive several preliminary facts.
The third section recalls the statements of the conjectures
and contains the proof of the main theorem.
The last section presents a short overview of known cases
of the Mumford--Tate conjecture and pointers to the literature.

\section{Preliminaries}

Let $k$ be a field of characteristic~$0$,
and fix a complex embedding $\sigma \colon k \into \CC$.
Our results do not depend on the choice of~$\sigma$.
We denote by $\Gamma_k$ the absolute Galois group $\Gal(\bar k/k)$,
where $\bar k$ is the algebraic closure of~$k$ in~$\CC$
along the embedding~$\sigma$.

If $C$ is a Tannakian category, over a field~$Q$
and $S$ is an object (or a collection of objects) of~$C$,
then we denote by $\langle S \rangle \subset C$ the smallest
Tannakian subcategory of~$C$ that contains~$S$. In other words, $\langle S
\rangle$ is the full subcategory of~$C$ that is the closure of~$S$ under
direct sums, tensor products, duals, and subquotients.

If $\omega \colon C \to \Vect{Q}$ is a fibre functor,
then we denote by $\group{G}\omega(S)$ the (pro)-algebraic group scheme
$\Aut^\otimes(\omega|_{\langle S \rangle})$ over~$Q$.
If $S = C$, we simply write $\group{G}\omega$ for $\group{G}\omega(C)$.
If $\langle S' \rangle \subset \langle S \rangle$,
then we get a natural surjection $\group{G}\omega(S) \to \group{G}\omega(S')$.

\medskip\noindent
\textbf{The category of motives.}\quad
Let $\Mot{k}$ denote the category of motives
in the sense of Yves Andr\'e~\cite{MR1423019}.
(To be precise,
we use the category of smooth projective $k$-schemes as ``base pieces'',
and we use singular cohomology relative to~$\sigma$ as ``reference cohomology''.
See \S2.1 of~\cite{MR1423019}.)
Alternatively,
one could use the theory of motives for absolute Hodge cycles~\cite{MR0654325};
this would not alter any of the following statements or proofs.
Recall from th\'eor\`eme~0.4 of~\cite{MR1423019} that $\Mot{k}$ is a
graded, polarisable, semisimple Tannakian category over~$\QQ$.

The complex embedding $\sigma \colon k \into \CC$
induces a realisation functor $r_\sigma \colon \Mot{k} \to \HS{\QQ}$
to the category of $\QQ$-Hodge structures.
Every prime number~$\ell$
induces a realisation functor $r_\ell \colon \Mot{k} \to \Rep{\QQl}(\Gamma_k)$
to the category of $\ell$-adic representations of~$\Gamma_k$.
We will denote by
$u_\sigma \colon \HS{\QQ} \to \Vect{\QQ}$
and $u_\ell \colon \Rep{\QQl}(\Gamma_k) \to \Vect{\QQl}$
the respective forgetful functors.
By Artin's comparison theorem (expos\'e~\textsc{ix} of~\cite{SGA4-3}),
we obtain a natural isomorphism of functors
$(- \otimes \QQl) \circ u_\sigma \circ r_\sigma \cong u_\ell \circ r_\ell$.

\medskip\noindent
\textbf{Artin motives and abelian motives.}\quad
A motive is called an \emph{Artin} motive if it is isomorphic to an object of
the Tannakian subcategory of~$\Mot{k}$ generated by the motives $\HH(X)$,
where $X$ ranges over all finite \'etale $k$-schemes.
(See example~(ii) after~\S4.5 of~\cite{MR1423019}.)

A motive is called an \emph{abelian} motive if it is isomorphic to an object of
the Tannakian subcategory of~$\Mot{k}$
generated by Artin motives and the motives $\HH(X)$,
where $X$ ranges over all abelian varietes over~$k$.
(See \S6.1 of~\cite{MR1423019}.)

\medskip\noindent
\textbf{Tate triples.}
\quad
Recall from~\S5 of~\cite{Del-Mil}, that a \emph{Tate triple} over a field~$Q$
is a triple $(C, w, T)$ consisting of
a Tannakian category~$C$ over~$Q$,
a cocharacter $w \colon \Gm \to \Aut^\otimes(\id_C)$, and
an invertible object~$T$ of weight~$-2$ (called the \emph{Tate object}).
If $\omega \colon C \to \Vect{Q}$ is a fibre functor,
then we obtain a natural surjection
$t \colon \group{G}\omega \to \group{G}\omega(T) = \Gm$
such that the composition $t \circ w = -2 \in \Hom(\Gm, \Gm)=\ZZ$.
Following \S5 and~\S13 of~\cite{Serre94},
we denote the kernel of~$t$ with a subscript~$1$, as in $\sp\group{G}\omega$.
(Note that~\cite{Del-Mil} uses the notation $G_0$ instead.)

As before, let $S$ be an object (or a collection of objects) of~$C$.
The image of $\sp\group{G}\omega$ in $\group{G}\omega(S)$
is denoted with $\sp\group{G}\omega(S)$.
We will now give an alternative description of $\sp\group{G}\omega(S)$.
We order $\NN = \ZZ_{\ge 0}$ by the divisibility relation,
so that $0$ is a maximal element.
Let $m$ be the smallest natural number for this ordering
such that $\langle S \rangle$ contains an object isomorphic to $T^{\otimes m}$.
The quotient $\group{G}\omega(S)/\sp\group{G}\omega(S)$ naturally
receives a surjective homomorphism from $\group{G}\omega(T) = \Gm$,
which factors through $\group{G}\omega(T^{\otimes m})$ by minimality of~$m$.
We thus obtain the following diagram with exact rows:
\[
 \begin{tikzcd}
  0 \rar & \sp\group{G}\omega \rar[hook] \dar[dotted] &
  \group{G}\omega \rar[two heads] \dar[two heads] &
  \group{G}\omega(T) = \Gm \dar[two heads,"(\_)^m"] \rar & 0\\
  0 \rar & K \rar[hook] \dar[dotted,hook] &
  \group{G}\omega(S) \rar[two heads]
  \arrow[transform canvas={xshift=0.3ex},-]{d} \arrow[transform canvas={xshift=-0.4ex},-]{d} &
  \group{G}\omega(T^{\otimes m}) \dar[two heads] \rar & 0 \\
  0 \rar & \sp\group{G}\omega(S) \rar[hook] &
  \group{G}\omega(S) \rar[two heads] &
  \group{G}\omega(S)/\sp\group{G}\omega(S) \rar & 0
 \end{tikzcd}
\]
By a formal argument, we see that the natural map
$\sp\group{G}\omega \onto \sp\group{G}\omega(S)$
factors through~$K$ as indicated by the dotted arrows.
We conclude that $\sp\group{G}\omega(S)$ is the kernel of the natural map
$\group{G}\omega(S) \to \group{G}\omega(T^{\otimes m})$.

The categories $\Mot{k}$, $\HS{\QQ}$,
and $\Im(r_\ell) \subset \Rep{\QQl}(\Gamma_k)$
are all equipped with the structure of a Tate triple
that is compatible with the realisation functors $r_\sigma$ and $r_\ell$.

\medskip\noindent
\textbf{Tannaka groups associated with motives.}\quad
We now specialise the previous discussion to the Tate triples
$\Mot{k}$, $\HS{\QQ}$ and $\Im(r_\ell) \subset \Rep{\QQl}(\Gamma_k)$.
(We restrict to the essential image $\Im(r_\ell)$ in order to obtain
a natural grading by weights.)

In the case of $\Mot{k}$,
let $\omega$ be the fibre functor
$u_\sigma \circ r_\sigma \colon \Mot{k} \to \Vect{\QQ}$.
We write $\mG{k}$ for $\group{G}\omega = \Aut^\otimes(\omega)$.
The notations $\mG{k}(M)$, $\sp\mG{k}$, and $\sp\mG{k}(M)$
are analogously defined.

For the Hodge realisation we use the following notation:
$\rGs = \group{G}{u_\sigma} = \Aut^\otimes(u_\sigma)$.
If $M$ is a motive, then we denote by $\rGs(M)$
the group scheme $\Aut^\otimes(u_\sigma|_{\langle r_\sigma(M) \rangle})$;
it is the \emph{Mumford--Tate group} of the Hodge structure~$r_\sigma(M)$.
The notations $\sp\rGs$ and $\sp\rGs(M)$ are similarly defined.
The functor $r_\sigma$ induces a morphism $\rGs \to \mG{k}$,
and $\rGs(M)$ may be identified with the image of $\rGs$ in~$\mG{k}(M)$.

For the $\ell$-adic realisation we use the following notation:
$\rGl = \group{G}{u_\ell}$, $\rGl(M)$, $\sp\rGl$, and $\sp\rGl(M)$.
The group scheme $\rGl(M)$ is
the so-called \emph{$\ell$-adic monodromy group} of~$M$.
It is the Zariski closure of the image of $\Gamma_k$ in $\GL(r_\ell(M))$.
Artin's comparison isomorphism induces a morphism $\rGl \to \mG{k,\QQl}$,
and $\rGl(M)$ may be identified with the image of $\rGl$ in~$\mG{k}(M)_{\QQl}$.

For any algebraic group~$\rG{}$,
we denote by $\pi_0\rG{}$ the component group of~$\rG{}$.
We write $\cc\mG{k}(M)$ (resp.~$\rGlc(M)$)
for the identity component of $\mG{k}(M)$ (resp.~$\rGl(M)$).
Note that $\rGs(M)$ is always a connected algebraic group.

\section{Main result}

\def\conjline#1:#2.{
\medskip
\indent\hbox{
 \hbox to  5em {$#1_\ell(M)$:\hfil}
 \hbox to 12em {#2,\hfil}
 \hbox to  5em {$#1(M)$:\hfil}
 \hbox to  6.5em {$\forall \ell,\ #1_\ell(M)$\hfil}
}
\medskip}

\begin{conjectures-non}
 Assume that the field $k$ is finitely generated as field,
 and let $M$ be a motive over~$k$.
 We recall the following conjectures.
 \begin{enumerate}
  \item \label{motivic-Tate-conjecture}
   A motivic analogue of the Tate conjecture:

   \conjline
   \motTC: $\rGl(M) = \mG{k}(M)_{\QQl}$.

   \noindent
   (N.b.: The ``classical'' $\ell$-adic Tate conjecture for
   $\ell$-adic cohomology of a smooth projective variety~$X/k$
   does not formally imply $\motTC_\ell(\HH(X))$,
   nor is the converse implication a formal fact:
   The conjecture $\motTC_\ell(\HH(X))$ expresses the assertion
   that all Tate classes in all cohomology groups of all powers of~$X$
   are motivated cycle classes.)
  \item \label{motivic-Sato--Tate-conjecture}
   The following is called the
   \emph{motivic Sato--Tate conjecture}
   in conj.~10.7 of~\cite{MR3502937}:

   \conjline
   \MST: $\sp\rGl(M) = \sp\mG{k}(M)_{\QQl}$.

  \item \label{motivic-Mumford--Tate-conjecture}
   A motivic version of the \emph{Mumford--Tate conjecture}:

   \conjline
   \MTC: $\rGlc(M) = \rGs(M)_{\QQl}$.

  \item \label{algebraic-Sato--Tate-conjecture}
   The \emph{algebraic Sato--Tate conjecture}
   (conj.~5.1(a,b) of~\cite{MR3502937}):

   \medskip
   {\narrower\it\noindent
   For every prime~$\ell$,
   there exists a natural-in-$k$ reductive algebraic group
   $\AST_k(M) \subset \GL(\omega(M))$ over~$\QQ$
   and a natural-in-$k$ isomorphism of group schemes
   $\sp\rGl(M) \cong \AST_k(M)_{\QQl}$.
   \par}
 \end{enumerate}
\end{conjectures-non}
Note that $\MST(M)$ is a more precise version
of the algebraic Sato--Tate conjecture:
it predicts that $\AST_k(M) = \sp\mG{k}(M)$.

\begin{theorem-non}\label{main-theorem}
Let $k$ be a finitely generated field of characteristic~$0$,
 and fix a complex embedding $\sigma \colon k \into \CC$.
 Let $M$ be an abelian motive over~$k$, and let $\ell$ be a prime number.
The following assertions are equivalent.

 \smallskip
 \noindent
 \begin{minipage}{.3\textwidth}
  \begin{enumerate}[labelwidth=2.5em,itemindent=!,labelsep=1.5em]
   \item[({\it i}\/)] $\motTC(M)$
   \item[({\it i}\/)\rlap{$_\ell$}] $\motTC_\ell(M)$
  \end{enumerate}
 \end{minipage}
 \begin{minipage}{.3\textwidth}
  \begin{enumerate}[labelwidth=2.5em,itemindent=!,labelsep=1.5em]
   \item[({\it ii}\/)] $\MST(M)$
   \item[({\it ii}\/)\rlap{$_\ell$}] $\MST_\ell(M)$
  \end{enumerate}
 \end{minipage}
 \begin{minipage}{.3\textwidth}
  \begin{enumerate}[labelwidth=2.5em,itemindent=!,labelsep=1.5em]
   \item[({\it iii}\/)] $\MTC(M)$
   \item[({\it iii}\/)\rlap{$_\ell$}] $\MTC_\ell(M)$
  \end{enumerate}
 \end{minipage}

 \smallskip
 \noindent
 In particular, $\MTC(M)$ implies the algebraic Sato--Tate conjecture for~$M$.
\end{theorem-non}
\noindent
\emph{Proof.}
Let $M$ be an motive over~$k$, and let $\ell$ be a prime number.
We start by proving $\MST_\ell(M) \iff \motTC_\ell(M)$.
This part of the proof does not use that $M$ is an abelian motive.
Consider the following commutative diagram with exact rows:
\[
\begin{tikzcd}
 0 \rar & \sp\rGl(M) \rar \dar[hook] & \rGl(M) \rar \dar[hook]
 & \rGl(M)/\sp\rGl(M) \rar \dar & 0 \\
 0 \rar & \sp\mG{k}(M)_{\QQl} \rar & \mG{k}(M)_{\QQl} \rar
 & \mG{k}(M)_{\QQl}/\sp\mG{k}(M)_{\QQl} \rar & 0
\end{tikzcd}
\]
The vertical arrow on the right is surjective, since
if the Tannakian subcategory $\langle M \rangle \subset \Mot{k}$ generated by~$M$
contains an object isomorphic to~$\mathds{1}(n)$, for some $n \in \ZZ$,
then the Tannakian subcategory $\langle r_\ell(M) \rangle$
of $\Im(r_\ell) \subset \Rep{\QQl}(\Gamma_k)$
contains an object isomorphic to~$\QQl(n)$.

We now make the following two observations.
\begin{enumerate}
 \item If the vertical arrow on the left is an isomorphism,
  then the inclusion in the middle must be an isomorphism
  by the five lemma.
 \item If the vertical arrow in the middle is an isomorphism,
  then $\langle M \rangle_{\QQl}$ is equivalent to
  $\langle r_\ell(M) \rangle$ and therefore both categories
  contain the same tensor powers of the Tate object.
  Therefore the vertical arrow on the right is an isomorphism,
  and so is the arrow on the left.
\end{enumerate}

Let us now prove $\MTC_\ell(M) \iff \motTC_\ell(M)$.
Consider the following commutative diagram with exact rows:
\[
\begin{tikzcd}
 0 \rar & \cc\rGl(M) \rar \dar[hook] & \rGl(M) \rar \dar[hook]
 & \pi_0\rGl(M) \rar \dar & 0 \\
 0 \rar & \cc\mG{k}(M)_{\QQl} \rar & \mG{k}(M)_{\QQl} \rar
 & \pi_0\mG{k}(M)_{\QQl} \rar & 0
\end{tikzcd}
\]
Now we argue as follows:
\begin{enumerate}
 \item Observe that if the inclusion in the middle is an isomorphism,
  then the inclusion on the left is an isomorphism,
  by definition of the identity component.
 \item Claim: The arrow $\pi_0\rGl(M) \to \pi_0\mG{k}(M)_{\QQl}$ is surjective.
  Indeed, since $\pi_0\mG{k}(M)$ is a finite group,
  we may view it as the motivic Galois group of some Artin motive~$N$.
  The image of $\pi_0\rGl(M)$ in $\pi_0\mG{k}(M)_{\QQl} = \mG{k}(N)_{\QQl}$
  is exactly $\rGl(N)$.
  The equality $\rGl(N) = \mG{k}(N)_{\QQl}$ is well-known;
  see for example remark~6.18 on page~211 of~\cite{Del-Mil}.

 \item Once again, the five lemma shows that if the inclusion on the left
  is an isomorphism and the right arrow is a surjection,
  then the inclusion in the middle is an isomorphism.
 \item To finish the proof, we assume that $M$ is an abelian motive.
  Under this assumption
  the canonical inclusion $\rGs(M) \into \cc\mG{k}(M)$ is an isomorphism.
  Indeed, by th\'eor\`eme~0.6.2 of~\cite{MR1423019}
  we know that $\rGs(M) \cong \mG{\CC}(M_\CC)$.
  We also have an isomorphism $\mG{\CC}(M_\CC) \cong \mG{\bar k}(M_{\bar k})$,
  see th\'eor\`eme~0.6.1
  and remarque~(ii) after th\'eor\`eme~5.2 of~\cite{MR1423019}.
  Thus it remains to show $\mG{\bar k}(M_{\bar k}) = \cc\mG{k}(M)$.
  \looseness=-1

  By example~(ii) in~\S4.6 of~\cite{MR1423019}
  we know that the quotient $\mG{k}(M)/\mG{\bar k}(M_{\bar k})$
  is a quotient of $\Gamma_k$.
  Since it is also a quotient of an algebraic group of finite type,
  this quotient is a finite group,
  and therefore isomorphic to $\pi_0\mG{k}(M)$.
\end{enumerate}
We conclude that $\MTC_\ell(M) \iff \motTC_\ell(M)$.
Finally, for abelian motives
the Mumford--Tate conjecture is independant of~$\ell$,
by corollary~7.6 of~\cite{1706.09444v2}:
it proves the implication $\MTC_\ell(M) \implies \MTC(M)$.
Altogether, this proves the theorem.
\hfill\qed

\begin{remark-non}
 Note that for motives~$M$ of the form $M = \HH^1(X)$,
 where $X$ is an abelian variety over~$k$,
 it was already known that the Mumford--Tate conjecture is independent of~$\ell$
 by work of Larsen and Pink, see theorem~4.3 of~\cite{LP95}.

 Hence we only need to refer to corollary~7.6 of~\cite{1706.09444v2}
 in the final step of the proof
 to show independence for arbitrary abelian motives~$M$.
\end{remark-non}

\smallskip

\section{New instances of the algebraic Sato--Tate conjecture}

The main theorem of our paper asserts that for abelian varieties
the Mumford--Tate conjecture is equivalent to the motivic Sato--Tate conjecture
and therefore implies the algebraic Sato--Tate conjecture.
In this section we give a non-exhaustive presentation
of some known results on the Mumford--Tate conjecture for abelian varieties.
We refer the reader to the survey paper of Moonen
for further details~\cite{Moonen}.

For instance, one of the first examples
where the Mumford--Tate conjecture is true
is the case of abelian varieties of CM~type~\cite{Pohlmann}.
Serre proved that the Mumford--Tate conjecture
is true for elliptic curves~\cite{Serre72}.
Moreover generally, the conjecture is known to be true
for abelian varieties of dimension less than or equal to~$3$
and also
for simple abelian varieties of prime dimension~\cite{Tankeev83,Ribet1983,Chi}.
By the work of Moonen and Zarhin~\cite{MZ99},
the Mumford--Tate conjecture is known
for abelian varieties of dimension less or equal to~$5$
that do not have a isogeny factor of dimension~$4$
with trivial endomorphism algebra.

Further results are known if we impose some conditions
on the endomorphism algebra of the abelian variety~$A$.
Indeed, if the endomorphism algebra of~$A$ is trivial,
and the dimension of~$A$ is an odd number,
Serre proved the Mumford--Tate conjecture~\cite{Serre2003}.
Several generalizations of this result were done afterwards by Chi~\cite{Chi2}
for abelian varieties with larger endomorphism algebra.
More results in this direction were proven by Pink~\cite{Pink98}.

Banaszak, Gajda and Kraso\'n~\cite{BGKIetII,BGK}
proved the Mumford--Tate conjecture
for some classes of abelian varieties of type~I,~II and~III
in the sense of Albert's classification.
Hindry and Ratazzi~\cite{HR} proved new instances
of the Mumford--Tate conjecture
for certain classes of abelian varieties of type~I and~II\@.
In~\cite{Can2017} the first author extends those results
to a larger class of abelian varieties of type~III.

Ichikawa~\cite{Ich} (resp.\@ Lombardo~\cite{LomHl}) proved that
under suitable conditions the Hodge group (resp.\@ $\ell$-adic Hodge group)
of a product of abelian varieties is the product of the Hodge groups
(resp.\@ $\ell$-adic Hodge groups).
Lombardo~\cite{LomHl} used these results to prove the Mumford--Tate conjecture
for arbitrary products of abelian varieties of dimension~$\le 3$.
Inspired by these results,
the second author~\cite{Commelin} proved that if two arbitrary abelian varieties
satisfy the Mumford--Tate conjecture,
then their product also satisfies this conjecture.

\printbibliography

\pagebreak

\end{document}